\documentclass{article}

\usepackage{amsmath}
\usepackage{amssymb}
\usepackage{graphicx}
\usepackage{doublespace}

\evensidemargin \oddsidemargin

\newtheorem{theorem}{Theorem}[section]
\newtheorem{lemma}[theorem]{Lemma}

\newtheorem{definition}[theorem]{Definition}

\numberwithin{equation}{section}

\def\P{\mathcal P}
\def\O{\mathcal O}

\begin{document}
\thispagestyle{empty}
\setcounter{page}{1}


\begin{center}
{\large\bf Necessity of Numerical Smoothness}

\vskip.20in

Tong Sun   
\\[2mm]

Department of Mathematics and Statistics\\
Bowling Green State University\\
Bowling Green, OH 43403 \\

\end{center}

\date{today}

\begin{spacing}{1.15}

{\footnotesize \noindent {\bf Abstract.} Numerical solutions of
differential equations are usually not smooth functions. However,
they should resemble the smoothness of the corresponding real
solutions in one way or another. In \cite{rumsey} and \cite{sun}, a
kind of spacial smoothness indicators was defined and subsequently
applied on the {\it a posteriori} error analysis. Here we prove that
the boundedness of those smoothness indicators is actually a
necessary condition for a piecewise polynomial function to
approximate a smooth function with optimal convergence rate. This
should help in validating the error analysis in \cite{rumsey} and
\cite{sun}. Moreover, the result of this paper provides an efficient
practical method to detect the loss of convergence rate due to the lack of
numerical smoothness, hence it serves as a criterion for the qualities of many numerical schemes. }

\vspace{0.2in}

{\bf keywords.} Numerical smoothness, smoothness indicator,
necessary condition.


\section{Introduction}

In the literature of numerical solutions of partial differential
equations, the smoothness of numerical solutions has not been a
popular concept. Of course, numerical solutions obtained by using
finite difference, finite element and finite volume methods are
typically not smooth functions. Therefore, it seems on the surface
that there is no smoothness whatsoever. However, it is widely known
that a scheme for solving a time dependent problem needs to be
dissipative and/or total variation diminishing. As a matter of
fact, these are actually concerning the smoothness of numerical
solutions. Although nonsmooth solutions (shocks, interfaces, etc.) are sometimes of great interest,
PDE solutions are usually at least piecewise smooth. A numerical scheme must
be able to approximate smooth pieces of solutions. Being dissipative and total variation
diminishing is certainly necessary. However, as shown in this paper, there are actually
stronger necessary conditions
to be satisfied if a scheme is hoped to converge at its desired optimal rate.

In \cite{rumsey} and \cite{sun}, a kind of spatial smoothness
indicators is proposed. These indicators are subsequently verified to be
bounded in the numerical experiments, which means that the computed numerical solutions
are ``numerically smooth". Most importantly, these smoothness indicators are
applied to the local error analysis, playing the role of higher order derivatives.
For the scalar nonlinear conservation laws studied in \cite{rumsey}
and \cite{sun}, obtaining the local error estimates in terms of the numerical
smoothness made it possible to do error propagation
analysis by directly using the $L^1$-contraction between entropy
solutions. Consequently, {\it a posteriori} error estimates of
optimal convergence rates and linear growth were obtained for the
RK-DG scheme and the WENO scheme. The key advantage of the methodology of  \cite{rumsey}
and \cite{sun} is that the smoothness indicators serve as a bridge
for bypassing the difficulty of proving any global property of a scheme, caused by the nonlinearities.

Here, we present the spatial smoothness indicator used in \cite{sun}.
We prove that the boundedness of the smoothness indicator is
actually a necessary condition for the numerical
solution to converge to any smooth function, including the real
solution, at the optimal convergence rate. In fact, the boundedness of the smoothness indicators 
is our choice of a global property to deal with in \cite{rumsey} and \cite{sun}. Since the numerical solutions are
computed by the very complex DG, WENO and Runge-Kutta schemes, we are
certainly unable to give any {\it a priori} proof for the
boundedness of the smoothness indicators. 
The computation of the smoothness indicators is actually the way 
to bypass the difficult proof. The necessity result of this article confirms that it is reasonable to expect 
the boundedness of the smoothness indicators. The result not only further
supports the error analysis of  \cite{rumsey} and \cite{sun}, but also supports the general methodology developed over there.

Beyond validating the methodology of \cite{rumsey} and \cite{sun}, the necessity
result provides an extremely efficient criterion on numerical
solutions of differential equations. Namely, if the smoothness
indicator computed from a numerical solution is too large, then the
numerical solution is certainly not converging to any smooth
function at the desired optimal rate. In another word, whenever the
smoothness indicators seem too large, either there is some kind of
non-smoothness in the PDE solution being approximated, or there is
something wrong in the numerical scheme being used.

The result is formulated in a 1-D uniform partition for the simplicity. One
can generalize the proof to higher dimensions and non-uniform triangulations.
However, as the first result of its kind, the investigation of a simple 1-D result suffices to show the necessity of numerical
smoothness.

\section{The main results}

Let $\P$ be the space of polynomials of degree $p$ or less, and $\P_h$ be
the space of piecewise polynomials of degree $p$ or less on the uniform
partition $a=x_0 < x_1 < \cdots < x_N=b$, $h=x_{i+1}-x_i=(b-a)/N$.
Let $u^R$ be a piecewise polynomial in $\P_h$.
The reason for using the superscript $R$ is to indicate that $u^R$ might be a reconstructed
numerical solution, possibly from nodal values, cell averages, or other forms of numerical
solutions. Define a smoothness indicator for any $u^R \in \P_h$ as
$$
S^p = (\bar{M}, \bar{D})
$$
with
$$
\bar{M} = (\tilde{M}_0, \tilde{M}_1, \cdots, \tilde{M}_{N-1}),
\qquad \text{where} \quad \tilde{M_i} = ( M^0_{i}, M^1_{i}, \cdots,
M^p_{i} ),
$$
and
$$
\bar{D} = (\tilde{D}_1, \cdots, \tilde{D}_{N-1}), \qquad
\text{where} \quad \tilde{D}_i = ( D^0_{i}, D^1_{i}, \cdots, D^p_{i}
).
$$
Furthermore,
$$
M_{i}^k = \frac{d^k}{dx^k} u^R(x_i^+), \qquad L_{i}^k =
\frac{d^k}{dx^k} u^R(x_i^-), \qquad J^k_i =  M_{i}^k -  L_{i}^k
$$
and, as in \cite{sun},
$$
D_{i}^k = J_{i}^k / h^{p+1-k}.
$$
In this article, $W^{k}_{q}$ is the standard notation of the Sobolev space of the functions, where
the $k$-th derivative is $L^q$-integrable. $H^k = W^{k}_{2}$. See \cite{brenner}, Chapter 2.

\begin{definition}
A piecewise polynomial $u^R \in \P_h$ is numerically
$W^{p+1}_{\infty}$-smooth in the  partition of cell size $h$, if there is a
constant $M_{\infty}$ such that all the components of $S^p$ are
bounded by $M_{\infty}$; $u^R$ is numerically $H^{p+1}$-smooth in
the partition if there is a constant $M_2$ such that
$\sum_{i=1}^{N-1} h \, [(D^0_i)^2 + (D^1_i)^2 + \cdots + (D^p_i)^2]
\le M_2^2$; $u^R$ is numerically $W^{p+1}_{1}$-smooth in the partition if
there is a constant $M_1$ such that $\sum_{i=1}^{N-1} h \, [|D^0_i|
+ |D^1_i| + \cdots + |D^p_i| ] \le M_1$.
\end{definition}

It is obvious that, for the piecewise polynomial $u^R$ to approximate a smooth function
in the interior of all the cells in the domain, one expects
$|M^k_{i}| \le \O(1)$ for all $k$ and $i$. It is also obvious that,
if $u^R$ is an approximation of a smooth function with the optimal
convergence rate, one can expect $|D^0_{i}| \le \O(1)$ for all $i$.
In the case $k=p$, it is obvious that $|D^p_i| \le \O(1)$ implies that
the piecewise constant function $\frac{d^p}{dx^p} u^R(x)$ has
bounded variation. In the main theorem below, we show the necessity
of the boundedness of the other components of the smoothness indicator. Let's start with a simple lemma.

\begin{lemma}
$$
Q(D^0,D^1,\cdots,D^p) =  \min_{\hat{v} \in \P} \left(
    \left\|\hat{v}+\frac{1}{2} \sum_{k=0}^{p} \frac{D^k}{k!}\xi^k \right\|^2_{L^2 (-\frac{1}{2},0)}
+ \left\|\hat{v}- \frac{1}{2} \sum_{k=0}^{p} \frac{D^k}{k!}\xi^k \right\|^2_{L^2 (0,\frac{1}{2})}  \right)
$$
is a positive definite quadratic form.
\end{lemma}
\noindent {\bf Proof.}  In the process of calculating the
coefficients of the polynomial $\hat{v} = \sum_{k=0}^{p} V_k \xi^k$
to attain the minimum, each $V_k$ must be a linear combination of
$D^0,D^1,\cdots,D^p$, hence the minimum is a quadratic form
of $D^0,D^1,\cdots,D^p$. The quadratic form is positive
definite because, if any entry $D^k$ of $(D^0,D^1,\cdots,D^p)$ is non-zero, then we have either
$V_k+\frac{1}{2}\frac{D^k}{k!} \ne 0$ or $V_k-\frac{1}{2}\frac{D^k}{k!} \ne 0$.
Due to the linear independence of $\{1,\xi,\cdots,\xi^p\}$, one of
the two terms in the minimum must be positive.  \#

\begin{theorem} \label{T2}
Suppose that $u \in H^{p+1}(a,b)$, and $u^R$ is a
piecewise polynomial function described as above. Then, there is a
constant $C_2 > 0$,  independent of $h$, $u$ and $u^R$, such that
\begin{equation} \label{L2D}
\|u - u^R\|_{L^2(a,b)} \ge  h^{p+1} \left[ \sqrt{\sum_{i=1}^{N-1} h
\, Q(D^0_i,D^1_i,\cdots,D^p_i)} - C_2 |u|_{H^{p+1}(a,b)} \right].
\end{equation}
\end{theorem}

\noindent {\bf Proof.} Let $\P_c$ be the space of piecewise
polynomials of degree $p$ or less on the partition
$a<x_{\frac{1}{2}}<x_{\frac{3}{2}}<\cdots<x_{N-\frac{1}{2}}<b$,
where $x_{i-\frac{1}{2}} = (x_{i-1}+x_i)/2$. Since $u \in
H^{p+1}(a,b)$, there is a constant $C_2 > 0$, independent of $h$ and
$u$, and a piecewise polynomial $u^I \in \P_c$, such that
\begin{equation} \label{uI}
\|u-u^I\|_{L^2 (a,b)} \le C_2 h^{p+1} |u|_{H^{p+1}(a,b)}.
\end{equation}
Now, since $\|u^I-u^R\|_{L^2 (a,b)} \le \|u^I-u\|_{L^2 (a,b)} + \|u-u^R\|_{L^2 (a,b)}$, we have
\begin{eqnarray}
&& \|u-u^R\|_{L^2 (a,b)} \ge \|u^I-u^R\|_{L^2 (a,b)} - \|u-u^I\|_{L^2 (a,b)} \nonumber \\
&& \ge \sqrt{ \sum_{i=1}^{N-1}  \|u^I-u^R\|^2_{L^2
(x_{i-\frac{1}{2}},x_{i+\frac{1}{2}})} } - C_2 h^{p+1}
|u|_{H^{p+1}(a,b)}. \label{jp}
\end{eqnarray}
Let $Q^k_i = (M^k_i + L^k_i)/2$ and $w(x) = \sum_{k=0}^p
\frac{Q^k_i}{k!} (x-x_i)^k$. Then, $u^R-w = \frac{1}{2} \sum_{k=0}^p
\frac{J^k_i}{k!} (x-x_i)^k$ for $x \in \Delta_i^+ =
(x_i,x_{i+\frac{1}{2}})$, and $u^R-w = - \frac{1}{2} \sum_{k=0}^p
\frac{J^k_i}{k!} (x-x_i)^k$ for $x \in \Delta_i^- =
(x_{i-\frac{1}{2}},x_i)$.
\begin{eqnarray}
&& \|u^I-u^R\|^2_{L^2 (x_{i-\frac{1}{2}},x_{i+\frac{1}{2}})} \ge \min_{v \in \P}  \|v-u^R\|^2_{L^2 (x_{i-\frac{1}{2}},x_{i+\frac{1}{2}})} \nonumber \\
&& =  \min_{v \in \P}  \|v-(u^R-w)\|^2_{L^2 (x_{i-\frac{1}{2}},x_{i+\frac{1}{2}})} \nonumber \\
&& =  \min_{v \in \P} \left( \left\|v+\frac{1}{2} \sum_{k=0}^{p}
\frac{J^k_i}{k!}(x-x_i)^k \right\|^2_{L^2(\Delta_i^-)}
                                 + \left\|v-\frac{1}{2} \sum_{k=0}^{p} \frac{J^k_i}{k!}(x-x_i)^k \right\|^2_{L^2 (\Delta_i^+)}  \right) \nonumber \\
&& =  h^{2p+3}  \min_{\hat{v} \in \P} \left(
\left\|\hat{v}+\frac{1}{2} \sum_{k=0}^{p} \frac{D^k_i}{k!}\xi^k
\right\|^2_{L^2 (-\frac{1}{2},0)}
              + \left\|\hat{v}-\frac{1}{2} \sum_{k=0}^{p} \frac{D^k_i}{k!}\xi^k \right\|^2_{L^2 (0,\frac{1}{2})}  \right) \nonumber \\
&& = h^{2p+3} \, Q(D^0_i, D^1_i,\cdots,D^p_i). \label{jpL}
\end{eqnarray}
Plugging (\ref{jpL}) into (\ref{jp}), we have proven (\ref{L2D}). \#

\begin{theorem} \label{T1}
If $u \in W^{p+1}_{1}(a,b)$, and $u^R$ is as described previously,
then, there are positive constants $C_1$ and $C_{12}^p$, independent
of $h$, $u$ and $u^R$, such that
\begin{equation} \label{L1D}
\|u - u^R\|_{L^1(a,b)} \ge  h^{p+1} \left[ C_{12}^p \sum_{1\le i \le
N-1} h \,\sqrt{ Q(D^0_i,D^1_i,\cdots,D^p_i)} - C_1
|u|_{W^{p+1}_{1}(a,b)} \right].
\end{equation}
\end{theorem}
\noindent {\bf Proof.} Since $u \in W^{p+1,1}(a,b)$, there is a
piecewise polynomial $u^I \in \P_c$, such that
$$
\|u-u^I\|_{L^1 (a,b)} \le  C_1 h^{p+1} |u|_{W^{p+1}_{1}(a,b)},
$$
where $C_1 >0$ is a constant independent of $h$ and $u$. Obviously,
\begin{eqnarray}
&& \|u-u^R\|_{L^1 (a,b)} \ge \|u^I-u^R\|_{L^1 (a,b)} - \|u-u^I\|_{L^1 (a,b)} \nonumber \\
&& \ge  \sum_{1 \le i \le N-1} \|u^I-u^R\|_{L^1
(x_{i-\frac{1}{2}},x_{i+\frac{1}{2}})}  - C_1 h^{p+1}
|u|_{W^{p+1}_{1}(a,b).} \label{jp1}
\end{eqnarray}
By the standard scaling argument, and (\ref{jpL}) which remains valid for the current $u^I$, it is easy to
prove that there is a constant $C_{12}^p$ such that
$$
\|u^I-u^R\|_{L^1 (x_{i-\frac{1}{2}},x_{i+\frac{1}{2}})} \ge C_{12}^p
h^{\frac{1}{2}} \|u^I-u^R\|_{L^2
(x_{i-\frac{1}{2}},x_{i+\frac{1}{2}})} \ge C_{12}^p h^{p+2}
\sqrt{Q(D^0_i,D^1_i,\cdots,D^p_i)}.
$$
Plugging this into (\ref{jp1}), we have (\ref{L1D}) proven. \#


\begin{theorem} \label{Tinfinity}
If $u \in W^{p+1}_{\infty}[a,b]$, and $u^R$ is as described previously, then,
there is a constant $C_\infty>0$, independent of $h$, $u$ and $u^R$,
such that
\begin{equation} \label{L8D}
\|u - u^R\|_{L^{\infty}(a,b)} \ge  h^{p+1} \left[\max_{1\le i \le
N-1} \sqrt{ Q(D^0_i,D^1_i,\cdots,D^p_i)} - C_\infty
|u|_{W^{p+1}_{\infty}[a,b]} \right].
\end{equation}
\end{theorem}
\noindent {\bf Proof.}
Since $u \in W^{p+1}_{\infty}[a,b]$, there is a piecewise polynomial $u^I \in \P_c$, such that
$$
\|u-u^I\|_{L^{\infty} (a,b)} \le  C_\infty h^{p+1}
|u|_{W^{p+1}_{\infty}[a,b]},
$$
where $C_\infty >0$ is a constant independent of $h$ and $u$.
Obviously,
\begin{eqnarray}
&& \|u-u^R\|_{L^{\infty} (a,b)} \ge \|u^I-u^R\|_{L^{\infty} (a,b)} - \|u-u^I\|_{L^{\infty} (a,b)} \nonumber \\
&& \ge  \max_{1 \le i \le N-1}
\|u^I-u^R\|_{L^{\infty}(x_{i-\frac{1}{2}},x_{i+\frac{1}{2}})}  -
C_\infty h^{p+1} |u|_{W^{p+1}_{\infty}[a,b].} \label{jp8}
\end{eqnarray}
Now, let $U_i =  \|u^I-u^R\|_{L^{\infty}(x_{i-\frac{1}{2}},x_{i+\frac{1}{2}})}$, then, by using (\ref{jpL}) for the current $u^I$,
\begin{eqnarray}
&& U_i = \sqrt{ \int_{x_{i-\frac{1}{2}}}^{x_{i+\frac{1}{2}}} \frac{1}{h} U_i^2 dx } \ge  \sqrt{ \int_{x_{i-\frac{1}{2}}}^{x_{i+\frac{1}{2}}} \frac{1}{h} (u^I-u^R)^2 dx } \nonumber \\
&& =  \sqrt{ \frac{1}{h} \|u^I-u^R\|^2_{L^2 (x_{i-\frac{1}{2}},x_{i+\frac{1}{2}})}  }
\ge h^{p+1} \sqrt{ Q(D^0_i,D^1_i,\cdots,D^p_i)}. \label{D8}
\end{eqnarray}
By plugging (\ref{D8}) into (\ref{jp8}), we have proven (\ref{L8D}).
\#

\section{Conclusion remarks}

According to Theorem \ref{T2}, in order to have
$$
\|u-u^R\|_{L^2(a,b)} \le \O (h^{p+1})
$$
for any function $u \in H^{p+1}(a,b)$, $u^R$ must be numerically $H^{p+1}$-smooth. That is
\begin{equation} \label{Hsmoothness}
\sum_{i=1}^{N=1} h\,  Q(D^0_i,D^1_i,\cdots,D^p_i) \le \O(1).
\end{equation}
According to Theorem \ref{T1}, in order to have
$$
\|u-u^R\|_{L^1(a,b)} \le \O (h^{p+1})
$$
for any function $u \in W^{p+1}_{1}(a,b)$, $u^R$ must be numerically
$W^{p+1}_{1}$-smooth. That is
\begin{equation} \label{Hsmoothness}
\sum_{i=1}^{N=1} h\, \sqrt{ Q(D^0_i,D^1_i,\cdots,D^p_i) } \le \O(1).
\end{equation}
According to Theorem \ref{Tinfinity}, in order to have
$$
\|u-u^R\|_{L^{\infty}(a,b)} \le \O (h^{p+1})
$$
for any function $u \in W^{p+1}_{\infty}[a,b]$, $u^R$ must be numerically $W^{p+1}_{\infty}$-smooth. That is
\begin{equation} \label{Csmoothness}
Q(D^0_i,D^1_i,\cdots,D^p_i) \le \O(1)
\end{equation}
for all $i$. Because $Q(d_0,d_1,\cdots,d_p)$ is a positive definite quadratic form, the last inequality is equivalent to that $|D^k_i| \le \O(1)$ for all $k$ and $i$. This is what we mean by the necessity of the numerical smoothness. To be more explicit, $|D^k_i| \le \O(1)$ is equivalent to $|J^k_i| \le \O(h^{p+1-k})$. That is, the jumps of the $k$-th derivative of $u^R$ need to be as small as $\O(h^{p+1-k})$, for $k=0,1,\cdots,p$.

Although the results and the proofs of the theorems seem to be very
simple, the author believes that the impact of the theorems could be
far reaching. For a numerical solution of a time-dependent PDE, the
theorem implies that a fully-discrete scheme must enforce numerical
smoothness, otherwise the convergence of optimal rate must have been
lost. Consequently, some of the traditional numerical stability
notions are possibly inadequate, if they do not enforce numerical
smoothness.

\bibliographystyle{amsplain}

\end{spacing}

\end{document}